\documentclass{amsart}
\usepackage{amsmath,amssymb}
\usepackage[dvips]{graphics}
\usepackage[all]{xy}

\makeatletter
\@addtoreset{figure}{section}
\def\@thmcountersep{-}
\makeatother

\numberwithin{equation}{section}

\begin{document}

\title{Unknotting numbers of diagrams of a given nontrivial knot are unbounded}

\author{Kouki Taniyama}
\address{Department of Mathematics, School of Education, Waseda University, Nishi-Waseda 1-6-1, Shinjuku-ku, Tokyo, 169-8050, Japan}
\email{taniyama@waseda.jp}


\date{}

\dedicatory{Dedicated to Professor Akio Kawauchi for his 60th birthday}


\begin{abstract}
We show that for any nontrivial knot $K$ and any natural number $n$ there is a diagram $D$ of $K$ such that the unknotting number of $D$ is greater than or equal to $n$. It is well known that twice the unknotting number of $K$ is less than or equal to the crossing number of $K$ minus one. We show that the equality holds only when $K$ is a $(2,p)$-torus knot.
\end{abstract}

\maketitle

\section{Introduction}	
Throughout this paper we work in the piecewise linear category. Let $L$ be a link in the 3-sphere ${\mathbb S}^3$ and $D$ a diagram of $L$ on the 2-sphere ${\mathbb S}^2$. It is well known that by changing over/under information at some crossings of $D$ we have a diagram of a trivial link. Let $u(D)$ be the minimal number of such crossing changes. Namely, there are some $u(D)$ crossings of $D$ such that changing them yields a trivial link diagram, and changing less than $u(D)$ crossings never yields a trivial link diagram. We call $u(D)$ the {\it unlinking number of $D$}. In the case that $D$ is a diagram of a knot $u(D)$ is called the {\it unknotting number of $D$}. The {\it unlinking number $u(L)$ of $L$} is defined by the minimum of $u(D)$ where $D$ varies over all diagrams of $L$. Namely we have the following equality.
\[
u(L)={\rm min}\{u(D)\mid D\mbox{ is a diagram of }L\}.
\]

\noindent
For a knot $K$ $u(K)$ is called the {\it unknotting number of $K$}. Then it is natural to ask whether or not the set $\{u(D)\mid D\mbox{ is a diagram of }L\}$ is bounded above. In \cite{Nakanishi1} Nakanishi showed that an unknotting number one knot $6_2$ has an unknotting number two diagram. Then he showed the following theorem in \cite{Nakanishi2}.

\vskip 3mm 
\noindent{\bf Theorem 1.1 \cite{Nakanishi2}.} {\it Let $K$ be a nontrivial knot. Then $K$ has a diagram $D$ with $u(D)\geq 2$.}

\vskip 3mm

\noindent
In this paper, as an extension of Theorem 1.1, we show the following theorem.

\vskip 3mm 
\noindent{\bf Theorem 1.2.} {\it Let $L$ be a nontrivial link. Then for any natural number $n$ there exists a diagram $D$ of $L$ with $u(D)\geq n$. \\That is, the set $\{u(D)\mid D\mbox{ is a diagram of }L\}$ is unbounded above.}

\vskip 3mm 

We note that Theorem 1.2 is an immediate consequence of the following proposition.

\vskip 3mm 
\noindent{\bf Proposition 1.3.} {\it Let $L$ be a nontrivial link and $D$ a diagram of $L$. Then there exists a diagram $D'$ of $L$ with $u(D')=u(D)+2$.}

\vskip 3mm

Let $c(D)$ be the number of crossings in $D$. We call $c(D)$ the {\it crossing number of $D$}. Then the {\it crossing number $c(L)$ of $L$} is defined by the minimum of $c(D)$ where $D$ varies over all diagrams of $L$. It is natural to ask the relation between $u(D)$ and $c(D)$, or $u(L)$ and $c(L)$. For a diagram $D$ of a knot $K$ other than a trivial diagram the following inequality is well-known. See for example \cite{Ozawa}.
\[
u(K)\leq u(D)\leq \frac{c(D)-1}{2}.
\]
In particular this inequality holds for a minimal crossing diagram $D$ of $K$ where $c(D)=c(K)$. Thus for any nontrivial knot $K$ we have the following inequality.
\[
u(K)\leq\frac{c(K)-1}{2}.
\]
It is also well known that the equality holds for $(2,p)$-torus knots. Conversely we have the following theorem.

\vskip 3mm 
\noindent{\bf Theorem 1.4.} {\it (1) Let $D$ be a diagram of a knot that satisfies the equality
\[
u(D)=\frac{c(D)-1}{2}.
\]
Then $D$ is one of the diagrams illustrated in Figure 1.1. Namely $D$ is a reduced alternating diagram of some $(2,p)$-torus knot, or $D$ is a diagram with just one crossing.

(2) Let $K$ be a nontrivial knot that satisfies the equality
\[
u(K)=\frac{c(K)-1}{2}.
\]
Then $K$ is a $(2,p)$-torus knot for some odd number $p\neq\pm1$. Namely only 2-braid knots satisfy the equality.}

\begin{figure}[htbp]
\begin{center}
\scalebox{0.5}{\includegraphics*{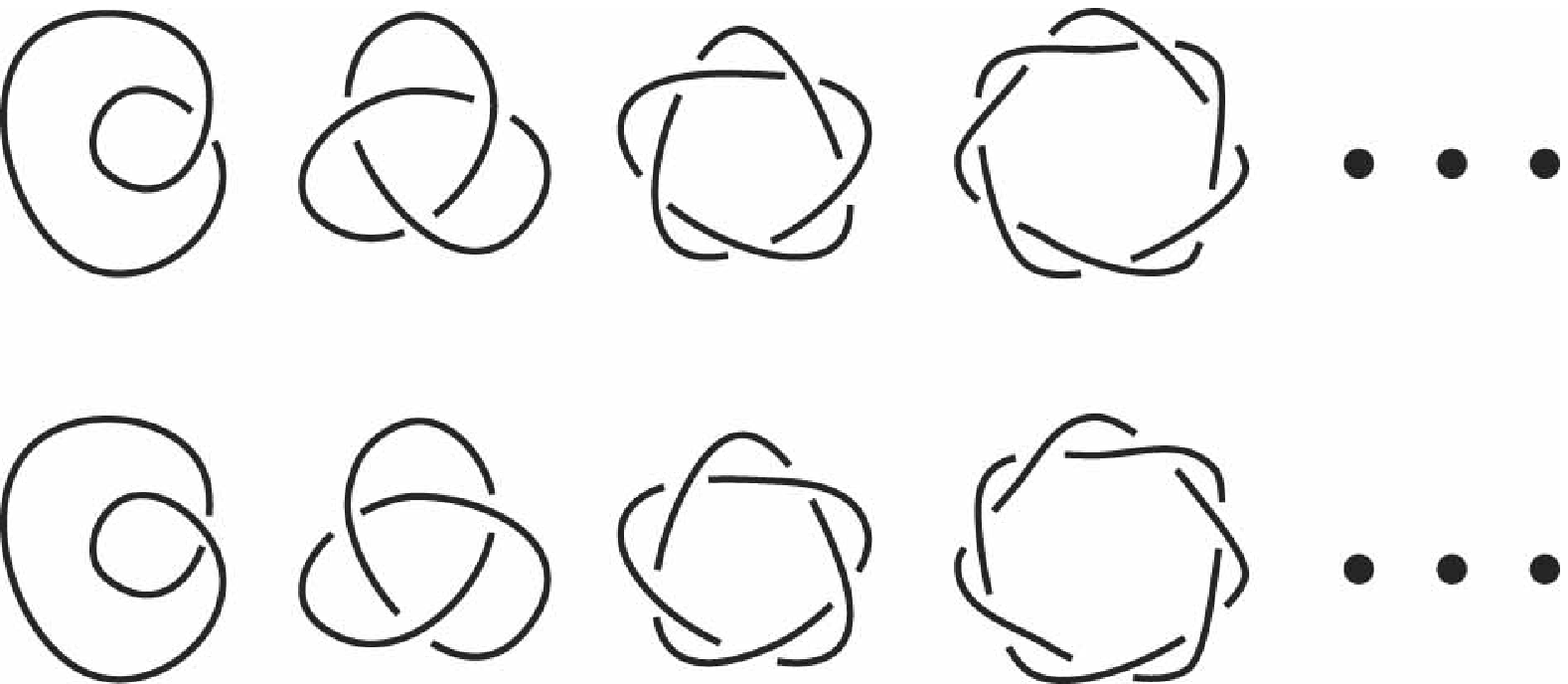}}
\end{center}
\centerline{Figure 1.1}
\end{figure}

\vskip 3mm

For links the situation is somewhat different. Let $D$  be a diagram of a link. Then the following inequality is well-known.
\[
u(L)\leq u(D)\leq \frac{c(D)}{2}.
\]
Thus for any link $L$ we have the following inequality.
\[
u(L)\leq\frac{c(L)}{2}.
\]
The following theorem shows that not only $(2,p)$-torus links but some other links satisfy the equality.

\vskip 3mm 
\noindent{\bf Theorem 1.5.} {\it (1) Let $D=\gamma_1\cup\cdots\cup\gamma_\mu$ be a diagram of a $\mu$-component link that satisfies the equality
\[
u(D)=\frac{c(D)}{2}.
\]
Then each $\gamma_i$ is a simple closed curve on ${\mathbb S}^2$ and for each pair $i,j$, the subdiagram $\gamma_i\cup\gamma_j$ is an alternating diagram or a diagram without crossings.

(2) Let $L$ be a $\mu$-component link that satisfies the equality
\[
u(L)=\frac{c(L)}{2}.
\]
Then $L$ has a diagram $D=\gamma_1\cup\cdots\cup\gamma_\mu$ such that each $\gamma_i$ is a simple closed curve on ${\mathbb S}^2$ and for each pair $i,j$, the subdiagram $\gamma_i\cup\gamma_j$ is an alternating diagram or a diagram without crossings.}

\vskip 3mm
Two examples of such links are illustrated in Figure 1.2. We note that for a link described in Theorem 1.5 the unlinking number equals the sum of the absolute values of all pairwise linking numbers. Let $a(L)$ be the ascending number of a link $L$ defined by Ozawa in \cite{Ozawa}. We also note here that if $K$ is a $(2,p)$-torus knot then $\displaystyle{a(K)=\frac{c(K)-1}{2}}$ and if $L$ is a link described in Theorem 1.5 then $\displaystyle{a(L)=\frac{c(L)}{2}}$. We do not know whether or not there exist other knots or links satisfying these equalities.

\begin{figure}[htbp]
\begin{center}
\scalebox{0.5}{\includegraphics*{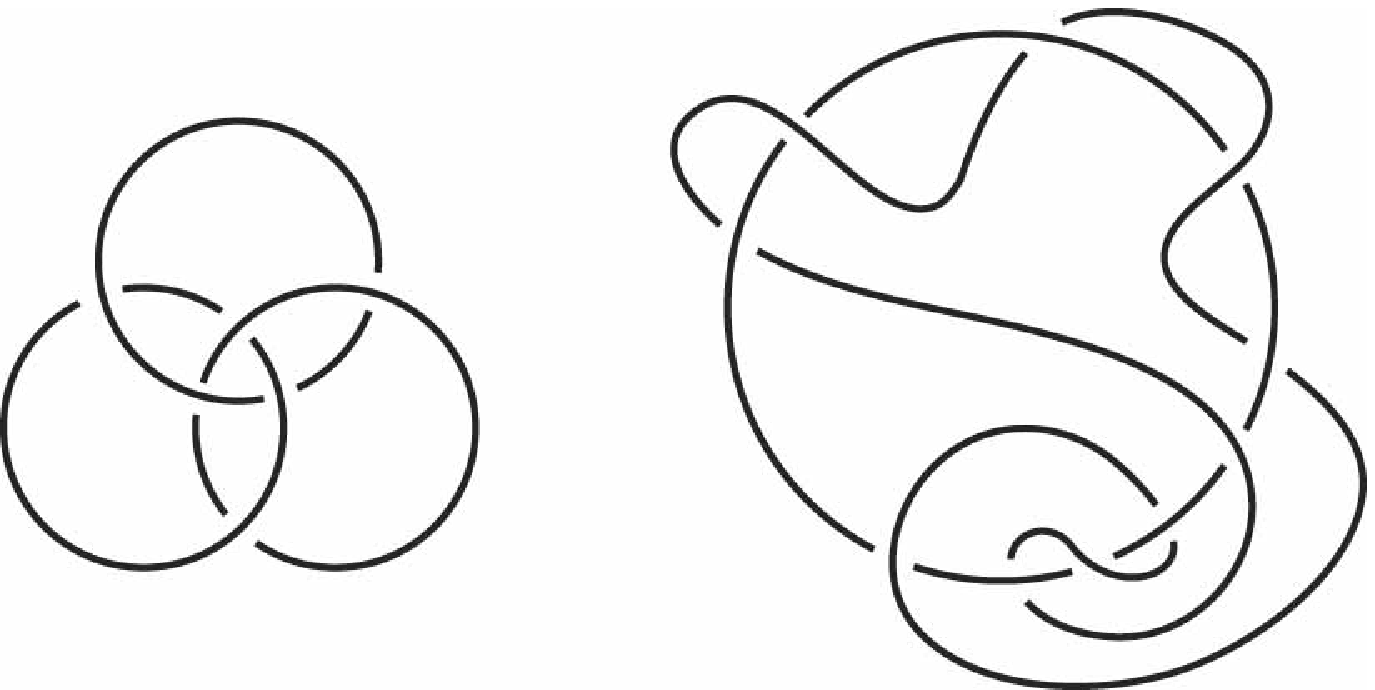}}
\end{center}
\centerline{Figure 1.2}
\end{figure}

\vskip 3mm

In section 2 we give a proof of Proposition 1.3. We then give some corollaries and questions. In section 3 we give proofs of Theorem 1.4 and Theorem 1.5.

\section{Proof of Proposition 1.3}
We first prepare the following two lemmas.

\vskip 3mm 
\noindent{\bf Lemma 2.1.} {\it Let $L$ be a link in ${\mathbb S}^3$ and $J$ a component of $L$. Let $L'$ be a link obtained from $L$ by adding some local knots to some components of $L$. Let $J'$ be the component of $L'$ that corresponds to $J$.  Then $L-J$ and $J$ are separable if and only if $L'-J'$ and $J'$ are separable.}

\vskip 3mm

\noindent{\bf Proof.} It is clear that if $L-J$ and $J$ are separable then $L'-J'$ and $J'$ are separable. We will show the converse. Suppose that $L'-J'$ and $J'$ are separable. Let $S$ be a separating sphere of them. Let $F_1,\cdots,F_k$ be decomposing spheres of $L'$. Namely each $F_i$ is a sphere intersecting $L'$ transversally at two points such that $F_i$ bounds a knotted ball-arc pair and if we replace each knotted arc by an unknotted arc then we have $L$. We may suppose that the intersection of each $F_i$ and $S$ are finitely many simple closed curves. Along an innermost disk on $F_i$ we cut $S$ into two spheres. We continue this until each $F_i$ has no intersection with spheres. Then we have a situation that there are some spheres, say $S_1,\cdots,S_l$ in ${\mathbb S}^3$ that are disjoint from $L'$ and each $F_i$. By considering black/white coloring of ${\mathbb S}^3$ by the spheres that is preserved under cutting operation we have that the component $N$ of ${\mathbb S}^3-(S_1\cup\cdots\cup S_l)$ containing $J'$ contains no other components of $L'$. The boundary of the closure of $N$ is a union of some of $S_1,\cdots,S_l$. After throwing away unnecessary spheres we pipe them and get a new separating sphere $S'$ of $L'-J'$ and $J'$ that is disjoint from each $F_i$. Then we have that $S'$ is also a separating sphere of $L-J$ and $J$. $\Box$

\vskip 3mm
In the following figures a right circle and a dotted line (resp. two dotted lines) inside it represents some 1-string (resp. 2-string) tangle possibly with some closed components.

\vskip 3mm 
\noindent{\bf Lemma 2.2.} {\it Let $L=J\cup M$ be a $\mu$-component link in ${\mathbb S}^3$ with $\mu\geq2$ as illustrated in {\rm Figure 2.1 (a)} where $J$ is a trivial knot, $M$ is a $\mu-1$-component link, $B$ is a 3-ball and the pair $(B,B\cap M)$ is a 2-string tangle with $\mu-2$ closed components. Suppose that $J$ and $M$ are separable. Then we have that the tangle $(B,B\cap M)$ is ambient isotopic relative to $\partial B$ to a tangle as illustrated in {\rm Figure 2.1 (b)} where $t_1$ and $t_2$ represents some 1-string sub-tangles possibly with some closed components.}

\begin{figure}[htbp]
\begin{center}
\scalebox{0.5}{\includegraphics*{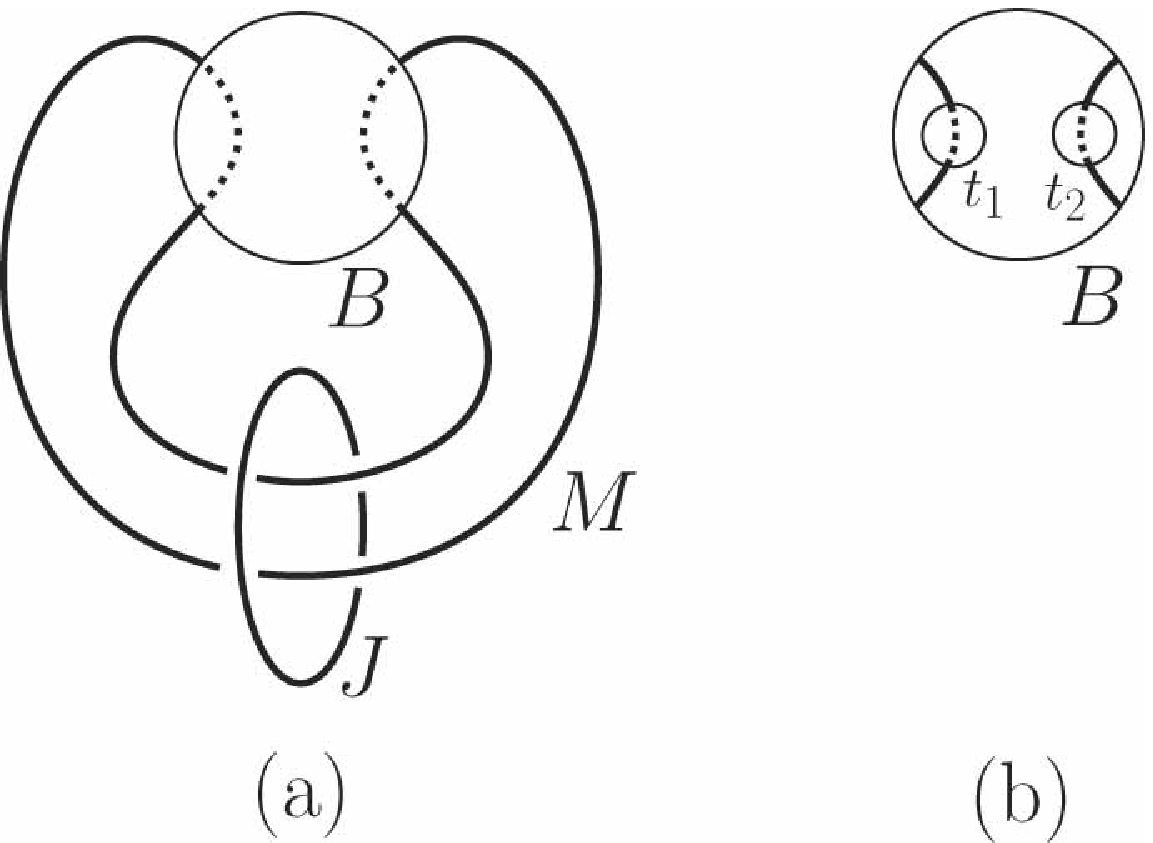}}
\end{center}
\centerline{Figure 2.1}
\end{figure}

\vskip 3mm

\noindent{\bf Proof.} Let $S=\partial B$ be the boundary of $B$. Suppose that $J$ and $M$ are separable. Then $J$ bounds a disk $\Delta$ that is disjoint from $M$. Since $J$ is homologically nontrivial in ${\mathbb S}^3-(B\cup M)$ we have that $\Delta$ cannot be disjoint from $B$. We may suppose that $\Delta$ intersects $S$ transversally and $\Delta\cap S$ is a disjoint union of finitely many simple loops, say $l_1,\cdots,l_k$. Let $d$ be an innermost disk in $\Delta$. Namely $d\cap(l_1\cup\cdots\cup l_k)=\partial d=l_i$ for some $i$. Let $B'$ be another 3-ball bounded by $S$. First suppose that $d$ is contained in $B'$. Since the tangle $(B',B'\cap L)$ is not a split tangle we have that the disk $d$ can be swept out of $B'$ by an ambient isotopy. Therefore we may suppose that $d$ is contained in $B$. If $d$ does not separates the strings then we can replace $\Delta$ by the disk with fewer intersection. Therefore we may suppose that $d$ separates the strings. Thus we have that the tangle $(B,B\cap M)$ is a split tangle. Therefore we have that $(B,B\cap M)$ is obtained from a rational tangle of some slope, say $q$ by adding local knots and closed components. Since $J$ and $M$ are separable we have that the lift of $M$ to the universal covering space of ${\mathbb S}^3-J$ which is homeomorphic to the 3-dimensional Euclidean space ${\mathbb R}^3$ is a splittable link of infinitely many components. Since it is as illustrated in Figure 2.2 (a) we have that a pair of adjacent components as illustrated in Figure 2.2 (b) is splittable. Then by Lemma 2.1 we have that the rational link of slope $1/q$ is splittable. Then by the classification of rational links \cite{Conway} \cite{Schubert} we have that this happens only when $1/q=0$. Namely we have $q=\infty$. $\Box$

\begin{figure}[htbp]
\begin{center}
\scalebox{0.5}{\includegraphics*{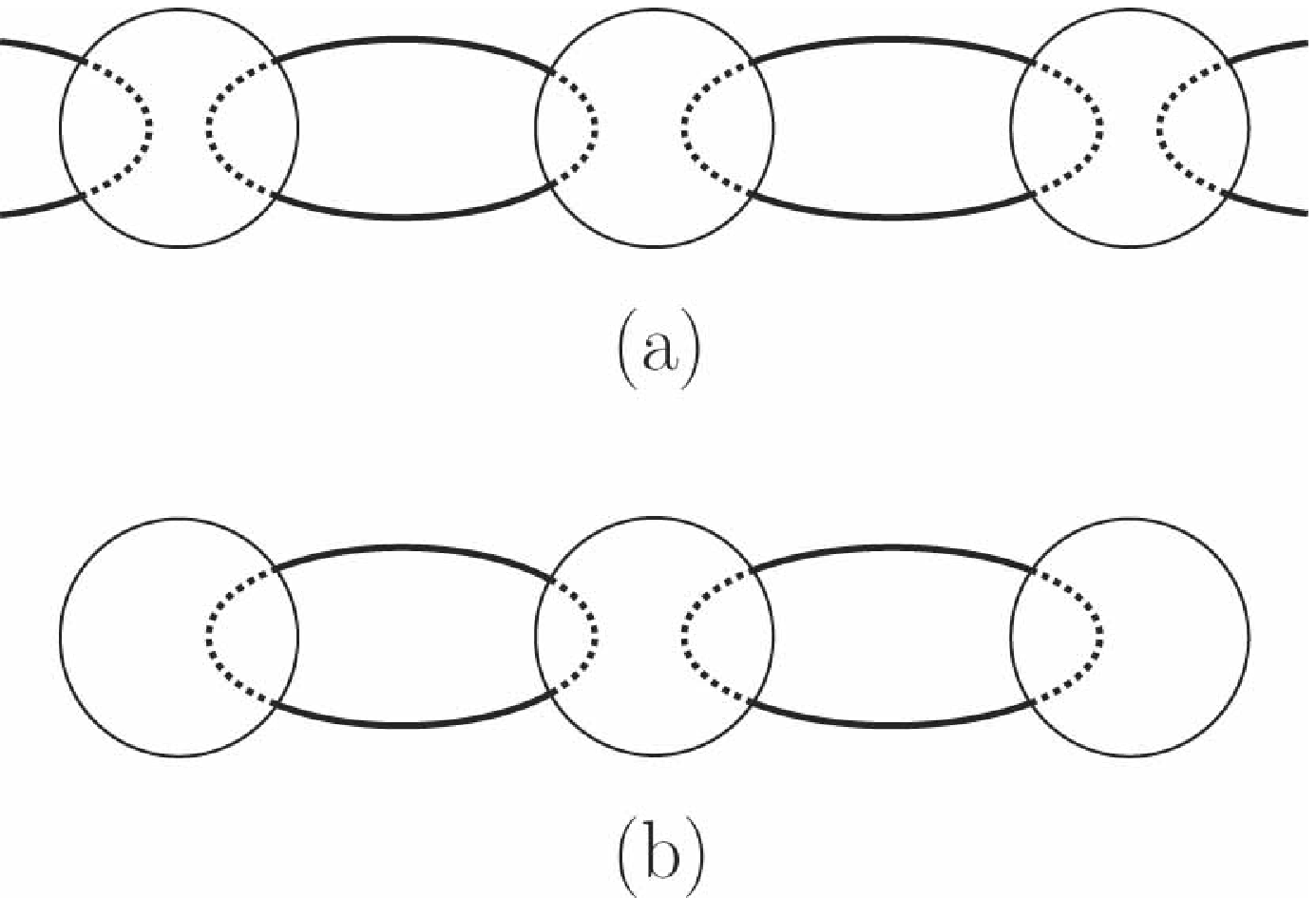}}
\end{center}
\centerline{Figure 2.2}
\end{figure}

\vskip 3mm

\noindent{\bf Proof of Proposition 1.3.} Let $C_1,\cdots,C_k$ be the set of all crossings of $D$. Set $l=u(D)$. By changing the order if necessary we may suppose that changing the crossings $C_1,\cdots,C_l$ yields a trivial link. Let $D'$ be a diagram obtained from $D$ by changing a neighbourhood of $C_i$ as illustrated in Figure 2.3 for each $i$ with $l\leq i\leq k$. It is clear that $D'$ is also a diagram of $L$. We note here that a prototype of this deformation is used in the proof of Theorem 1.1 in \cite{Nakanishi2}. We will show that $u(D')=u(D)+2$. By changing the crossings $C_1,\cdots,C_l,A_{l,1}$ and $A_{l,2}$ of $D'$ we have a trivial link. Thus we have $u(D')\leq u(D)+2$. Now suppose that $X$ is a set of crossings of $D'$ containing exactly $u(D')$ crossings such that changing all of them yields a trivial link $U$. We will show that $X$ contains exactly $u(D)+2$ crossings. First we show that the following four cases cannot happen.

\begin{figure}[htbp]
\begin{center}
\scalebox{1.1}{\includegraphics*{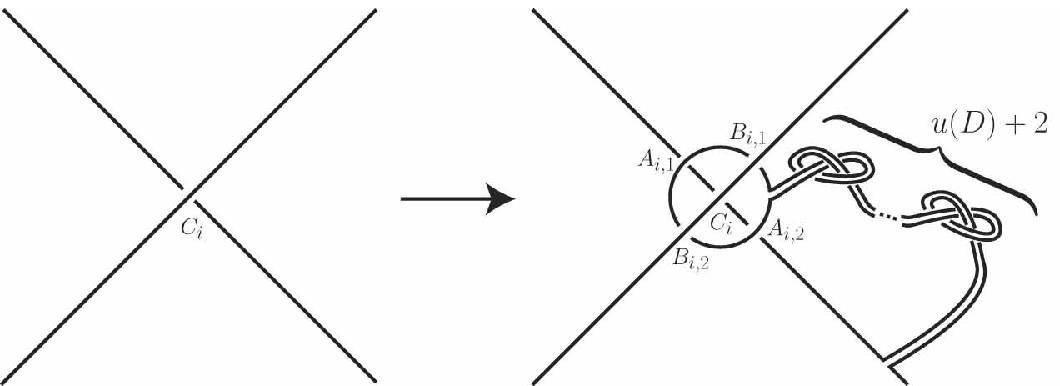}}
\end{center}
\centerline{Figure 2.3}
\end{figure} 

\vskip 3mm

Case 1. The number of elements of $\{A_{i,1},A_{i,2},B_{i,1},B_{i,2}\}\cap X$ is odd for some $i$.

Suppose for example that $C_i$ is a crossing of the same component of $L$ and $\{A_{i,1},A_{i,2},B_{i,1},B_{i,2}\}\cap X=\{B_{i,2}\}$. After changing the crossings in $X$ we have that at least one of $u(D)+2$ parallel trefoils in Figure 2.3 still alive. Then after an appropriate deformation that is fixed on a small neighbourhood of $B_{i,1}$ we have that $U$ is as illustrated in Figure 2.4 (a) or (b). We will show that the torus $T$ illustrated in Figure 2.4 is essential. Let $V$ be the solid torus in ${\mathbb S}^3$ bounded by $T$. If $T$ is inessential then a meridian of $T$ bounds a disk in $V$ that does not intersect $U$. Then we have that the lift of $U$ to the universal covering space of $V$ is a splittable link of infinitely many components. However it is easily seen that two adjacent components of them have linking number $1$ or $-1$. Thus we have that $T$ is essential. This contradicts the assumption that $U$ is a trivial link. Therefore this case cannot happen. If we take linking number into account in the case that $C_i$ is a crossing of some different components of $L$, all other possibilities of $\{A_{i,1},A_{i,2},B_{i,1},B_{i,2}\}\cap X$ can be checked in similar ways and we omit them.

\begin{figure}[htbp]
\begin{center}
\scalebox{0.5}{\includegraphics*{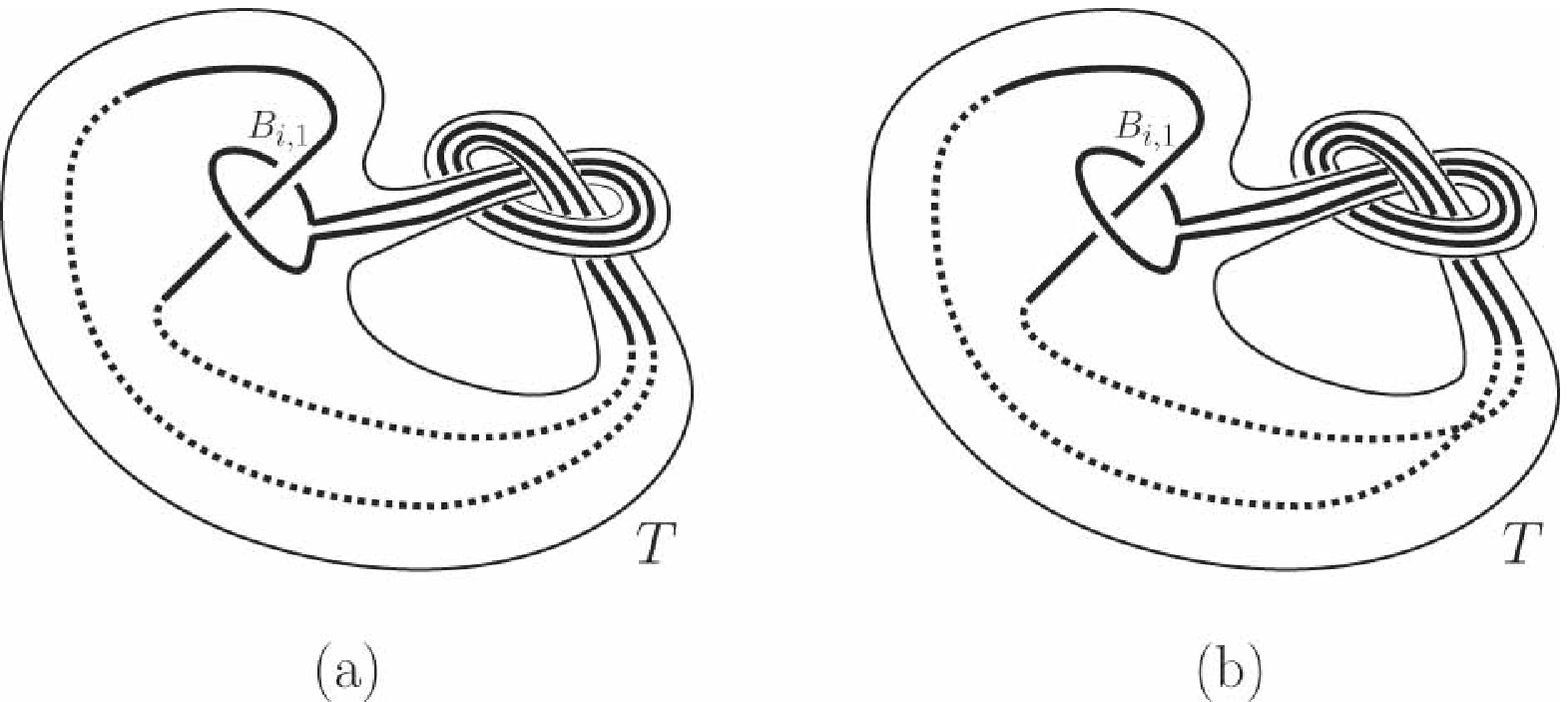}}
\end{center}
\centerline{Figure 2.4}
\end{figure} 

\vskip 3mm

Case 2. The set $\{A_{i,1},A_{i,2},B_{i,1},B_{i,2},C_i\}\cap X$ equals $\{C_i\}$ or \\$\{A_{i,1},A_{i,2},B_{i,1},B_{i,2}\}$ for some $i$.

Suppose for example that $C_i$ is a crossing of the same component of $L$ and $\{A_{i,1},A_{i,2},B_{i,1},B_{i,2},C_i\}\cap X=\{C_i\}$. Then we have that $U$ is as illustrated in Figure 2.5 (a) or (b). We only consider the case that $U$ is as illustrated in Figure 2.5 (a). The other case is essentially the same. Then it is deformed as illustrated in Figure 2.5 (c). We will show that the torus $T$ illustrated in Figure 2.5 (c) is essential. As in Case 1 we see the lift of $U$ to the universal covering space of the solid torus $V$ bounded by $T$. Then the adjacent components form a 2-component link $L'$ as illustrated in Figure 2.5 (d). It is sufficient to show that $L'$ is non-splittable. Then by Lemma 2.1 it is sufficient to show that the link $L''$ as illustrated in Figure 2.5 (e) is non-splittable. We may think that $L''$ is in ${\mathbb S}^3$ and we again consider the lift of $M$ to the universal covering space of ${\mathbb S}^3-J$. Then the adjacent components have linking number $\pm 1$. Therefore we have that the lift is not splittable. Then we have that $L''$ is non-splittable. The case that $C_i$ is a crossing of some different components of $L$ is similar and we omit it.

\begin{figure}[htbp]
\begin{center}
\scalebox{0.7}{\includegraphics*{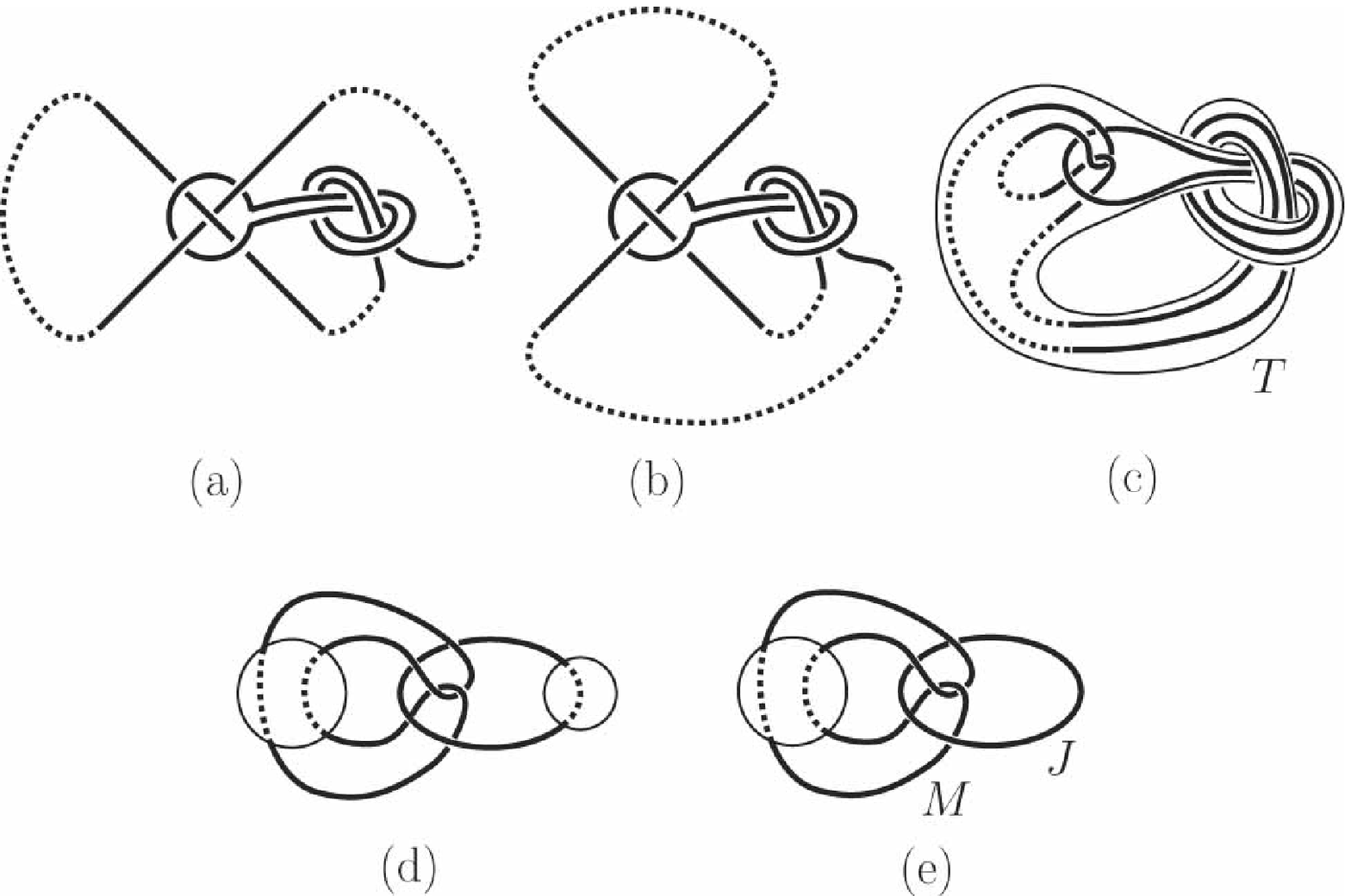}}
\end{center}
\centerline{Figure 2.5}
\end{figure} 

\vskip 3mm

Case 3. The set $\{A_{i,1},A_{i,2},B_{i,1},B_{i,2},C_i\}\cap X$ equals $\{A_{i,1},A_{i,2}\}$, $\{B_{i,1},B_{i,2}\}$ or \\
$\{A_{i,1},A_{i,2},B_{i,1},B_{i,2},C_i\}$ for some $i$.

In these cases we have the same link by 2 fewer crossing changes. This contradicts the assumption that $X$ contains exactly $u(D')$ crossings. Therefore these cases cannot happen.

\vskip 3mm

Case 4. Exactly one of $A_{i,1}$ and $A_{i,2}$ is contained in $X$ and exactly one of $B_{i,1}$ and $B_{i,2}$ is contained in $X$ for some $i$.

The case that $C_i$ is a crossing of some different components of $L$ is essentially the same as Case 1. Therefore we consider the case that $C_i$ is a crossing of the same component of $L$. Suppose for example that the set $\{A_{i,1},A_{i,2},B_{i,1},B_{i,2},C_i\}\cap X$ equals $\{A_{i,1},B_{i,1}\}$. Then we have that the link $U$ can be deformed into a form illustrated in Figure 2.6 (a) or (b). As in Case 1 we can show that the link illustrated in Figure 2.6 (a) has an essential torus $T$ by checking that an adjacent pair of components of a lift of the link to the universal covering space of the solid torus bounded by $T$ has linking number $\pm 2$. Now we consider the case that $U$ is as illustrated in Figure 2.6 (b). To see the situation clear we further deform $U$ as illustrated in Figure 2.6 (c). We will show that if the torus $T$ illustrated in Figure 2.6 (c) is not essential then the tangle $\tau$ illustrated in Figure 2.6 (c) is a tangle as illustrated in Figure 2.6 (d) up to ambient isotopy relative to the boundary of the 3-ball. Let $V$ be the solid torus bounded by $T$. Then the universal covering space of $V$ is as illustrated in Figure 2.6 (e). If $T$ is not essential then we have that a pair of adjacent lifts of $U$ in Figure 2.6 (e), where the left one is with the closed components contained in the tangle $\tau$, are separable. Then by Lemma 2.1 we have that the component outside the tangle $\tau$ of the link illustrated in Figure 2.6 (f) is separable from the rest. We note that the link illustrated in Figure 2.6 (f) is equivalent to the link described in Lemma 2.2. Then by Lemma 2.2 we have the desired conclusion. Then we have that $C_i$ is a nugatory crossing. Let $U'$ be the link obtained from $D'$ by changing all crossings in $X-\{A_{i,1},B_{i,1}\}$. Then we have that $U$ and $U'$ are ambient isotopic. This contradicts the assumption that $X$ contains exactly $u(D')$ crossings. Other cases of $\{A_{i,1},A_{i,2},B_{i,1},B_{i,2},C_i\}\cap X$ are quite similar and we omit them.

\begin{figure}[htbp]
\begin{center}
\scalebox{0.6}{\includegraphics*{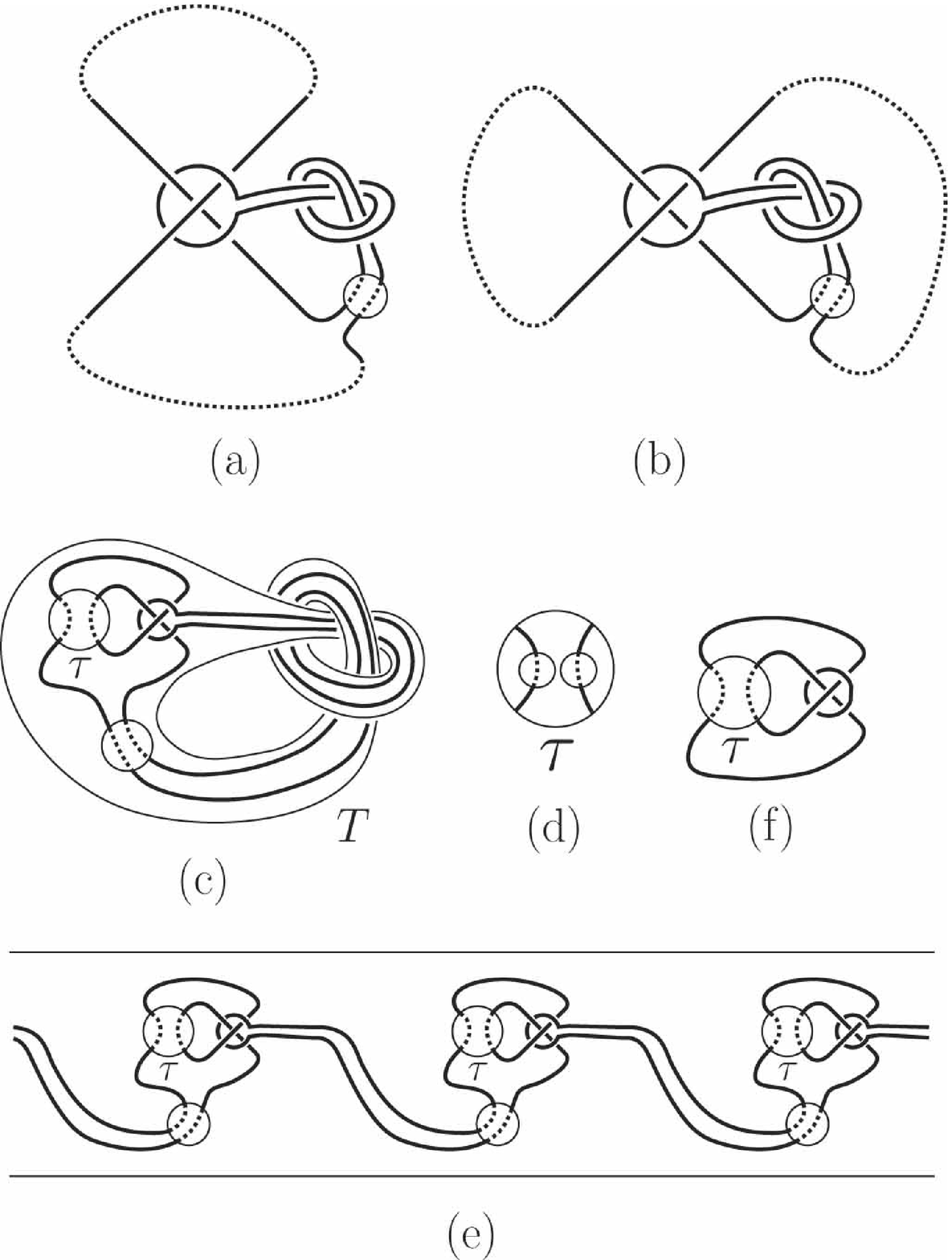}}
\end{center}
\centerline{Figure 2.6}
\end{figure} 

\vskip 3mm

Now we have that for each $i$ the set $\{A_{i,1},A_{i,2},B_{i,1},B_{i,2},C_i\}\cap X$ equals the empty set, $\{A_{i,1},A_{i,2},C_i\}$ or $\{B_{i,1},B_{i,2},C_i\}$. Then we have that the link $U$ is a composition of a link $W$ that is obtained from $D$ by changing all crossings in $\{C_1,\cdots,C_k\}\cap X$ and some local knots that arise by changing some crossings of parallel trefoils. Since $U$ is a trivial link we have that all these factors, in particular $W$, is trivial. Therefore we have that $D$ yields a trivial link by changing the crossings in $\{C_1,\cdots,C_k\}\cap X$. Thus we have that $\{C_1,\cdots,C_k\}\cap X$ contains at least $l=u(D)$ elements. Therefore at least one $C_i$ with $i\geq l$ is contained in $X$. Then $A_{i,1}$ and $A_{i,2}$, or $B_{i,1}$ and $B_{i,2}$ are also contained in $X$. Therefore we have that $X$ contains at least $u(D)+2$ crossings. Since we have shown $u(D')\leq u(D)+2$ we have that $X$ contains exactly $u(D)+2$ crossings. $\Box$

\vskip 3mm

As an immediate consequence of Proposition 1.3 we have the following corollary.

\vskip 3mm 
\noindent{\bf Corollary 2.3.} {\it Let $L$ be a nontrivial link. \\
Then the set $\{u(D)\mid D\mbox{ is a diagram of }L\}$ contains a set \\
$\{u(L)+2m\mid m\ \mbox{is a non-negative integer}\}$.}

\vskip 3mm 
\noindent{\bf Question 2.4.} {\it Let $L$ be a nontrivial link. Is the set $\{u(D)\mid D\mbox{ is a diagram of }L\}$ equals the set $\{u(L)+m\mid m\ \mbox{is a non-negative integer}\}$?}

\vskip 3mm 
The following proposition is a partial answer to Question 2.4.

\vskip 3mm 
\noindent{\bf Proposition 2.5.} {\it Let $L$ be an alternating link with $u(L)=1$. Suppose that $L$ has an alternating diagram $D_0$ with $u(D_0)=1$. \\
Then the set $\{u(D)\mid D\mbox{ is a diagram of }L\}$ equals the set of natural numbers \\
$\{u(L)+m\mid m\ \mbox{is a non-negative integer}\}$.}

\vskip 3mm 
\noindent{\bf Proof.} We may suppose without loss of generality that $D_0$ has no nugatory crossings. Let $D_1$ be a diagram of $L$ obtained from $D_0$ by replacing each neighbourhood of a crossing of $D_0$ as illustrated in Figure 2.7. It is easily seen that changing just one crossing of $D_1$ yields a link that has an alternating diagram without nugatory crossings. Therefore the link cannot be a trivial link. Thus we have that $u(D_1)\geq 2$. By changing two crossings as illustrated in Figure 2.7 we have the same effect of changing a crossing of the original diagram. This shows that $u(D_1)\leq 2$. Thus we have $u(D_1)=2$. Then by applying Proposition 1.3 to $D_0$ and $D_1$ repeatedly we have the desired conclusion. $\Box$

\begin{figure}[htbp]
\begin{center}
\scalebox{0.4}{\includegraphics*{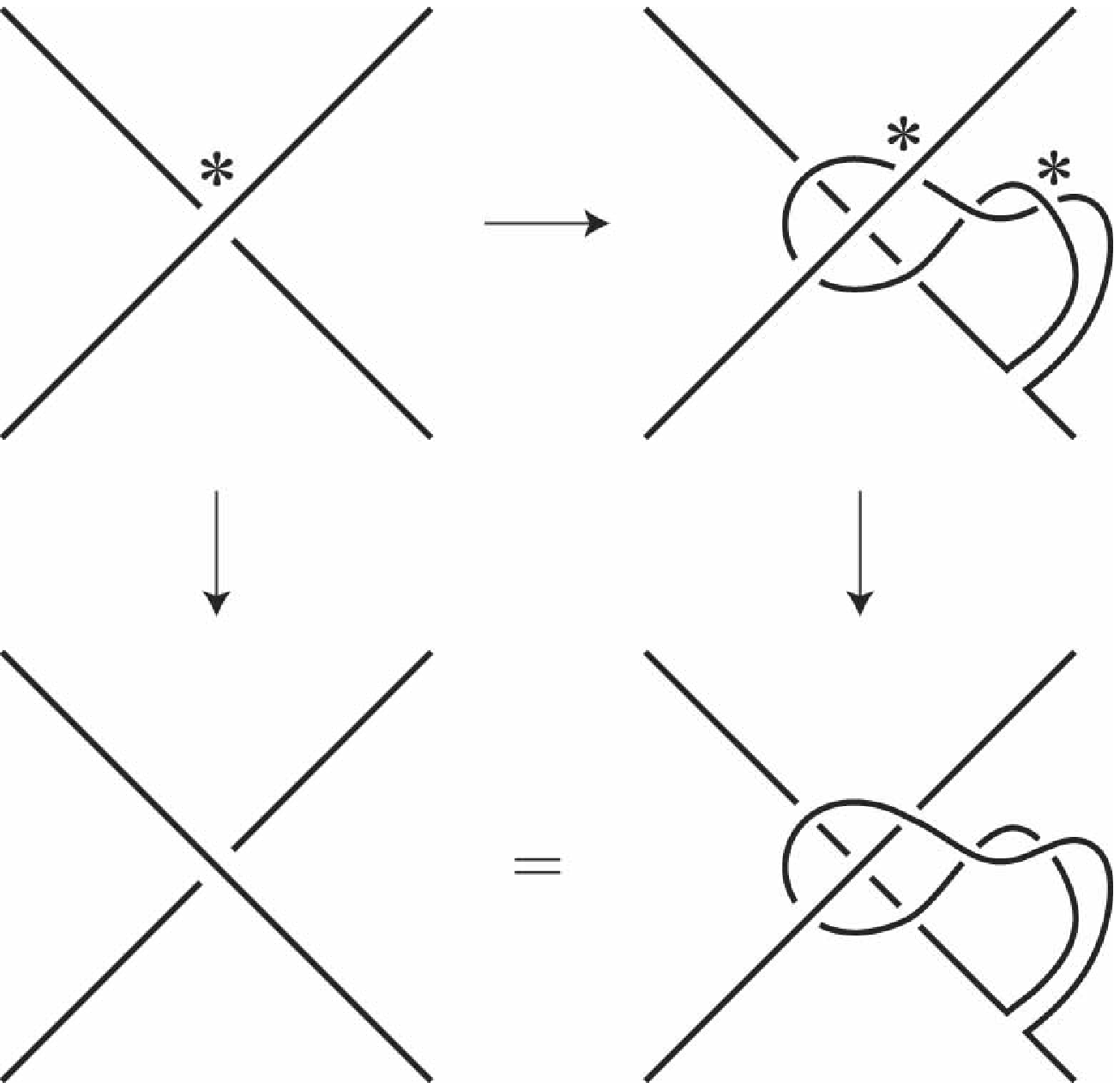}}
\end{center}
\centerline{Figure 2.7}
\end{figure} 

\vskip 3mm

\noindent{\bf Remark 2.6.} (1) Unknotting number one alternating knots with unknotting number one alternating diagrams are completely determined by Tsukamoto in \cite{Tsukamoto}. It is conjectured that every alternating diagram of an unknotting number one alternating knot has unknotting number one \cite{Nakanishi1} \cite{Kohn} \cite{Adams}. There is a more general conjecture by Bernhard \cite{Bernhard} and Jablan \cite{Jablan} that every nontrivial knot has a minimal crossing diagram such that changing a crossing in that diagram yields a knot with fewer unknotting number. See also \cite{Nakanishi1} \cite{Bleiler} \cite{Stoimenow} \cite{JS} \cite{DES} etc. for related problems.

(2) The construction illustrated in Figure 2.7 does not always work. A diagram of a trefoil knot and a diagram obtained from it by the construction of Figure 2.7 are illustrated in Figure 2.8. It is easy to check that both of them have unknotting number one.

(3) Our construction in the proof of Proposition 1.3 increases the crossing numbers of diagrams rapidly. Finding diagrams with fewer crossings will be an interesting problem.

\begin{figure}[htbp]
\begin{center}
\scalebox{0.7}{\includegraphics*{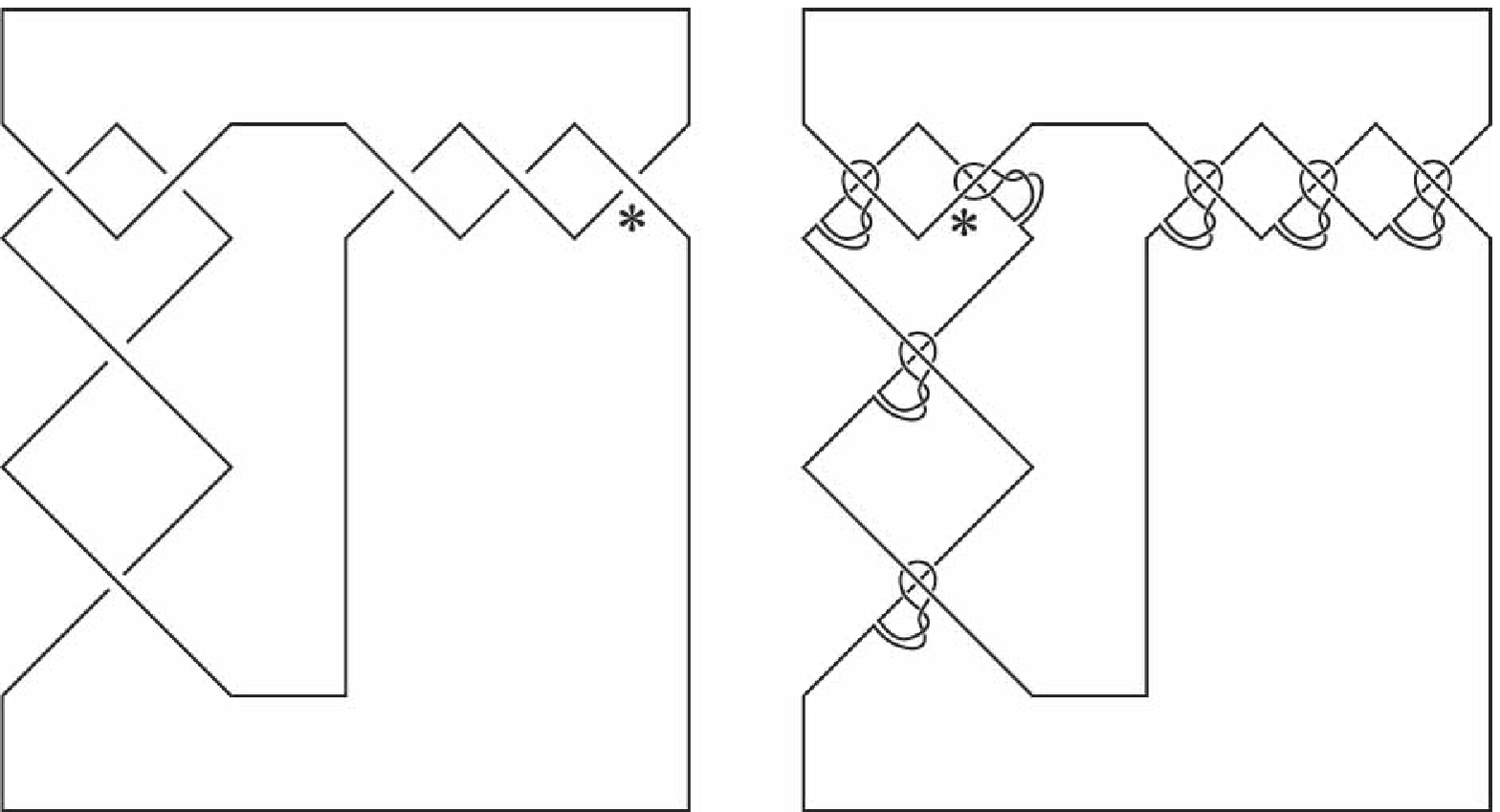}}
\end{center}
\centerline{Figure 2.8}
\end{figure} 

\vskip 3mm

\section{Proofs of Theorem 1.4 and Theorem 1.5}

\vskip 3mm 
\noindent{\bf Proof of Theorem 1.4.} (1) It is sufficient to show that if a knot diagram $D$ other than the trivial diagram is not in Figure 1.1 then $\displaystyle{u(D)<\frac{c(D)-1}{2}}$. Let $P$ be a crossing of $D$. By smoothing $D$ at $P$ we have a diagram $D'$ of a 2-component link. Let $\gamma_1$ and $\gamma_2$ be the components of $D'$. Suppose that one of them, say $\gamma_1$ is not a simple closed curve on ${\mathbb S}^2$. Let $Q$ be a self crossing of $\gamma_1$. By smoothing $D$ at $Q$ we have a diagram $D''$ of a 2-component link. Let $\gamma_3$ and $\gamma_4$ be the components of $D''$. We may suppose without loss of generality that $P$ is a self crossing of $\gamma_3$. See Figure 3.1. Then the diagram $\gamma_3$ is a nontrivial diagram of some knot. Therefore we have that $\displaystyle{u(\gamma_3)\leq\frac{c(\gamma_3)-1}{2}}$. If we change some crossings on $\gamma_4$ so that the part $\gamma_4$ is over other strings of $D$ and itself unknotted then we have a knot that has a diagram $\gamma_3$. Also we may change some crossings on $\gamma_4$ so that the part $\gamma_4$ is under other strings of $D$. Note that these two crossing changes are complementary on the crossings on $\gamma_4$. We choose one of them that have no more crossing changes than the other. Thus by changing no more than $\displaystyle{\frac{c(D)-c(\gamma_3)-1}{2}}$ crossings of $D$ we have a knot that has a diagram $\gamma_3$. Note that the key point here is that we do not need to change the crossing $Q$. Therefore we have that $\displaystyle{u(D)\leq u(\gamma_3)+\frac{c(D)-c(\gamma_3)-1}{2}}$. Therefore we have $\displaystyle{u(D)\leq\frac{c(D)-2}{2}}$. Thus we may suppose that both $\gamma_1$ and $\gamma_2$ are simple closed curves on ${\mathbb S}^2$. Now we trace $D$ on $\gamma_2$ starting from $P$ and see how it crosses with $\gamma_1$. If we find a situation as illustrated in Figure 3.2 then by replacing $P$ with $P'$ we have the previous situation. Then we finally have that the underlying projection of $D$ is the underlying projection of one of the diagrams illustrated in Figure 1.1. See Figure 3.3. Then it is clear that only alternating over/under crossing information satisfies the equality as desired.

(2) Let $K$ be a nontrivial knot that satisfies the equality $\displaystyle{u(K)=\frac{c(K)-1}{2}}$. Let $D$ be a minimal crossing diagram of $K$. Since $u(K)\leq u(D)$ and $\displaystyle{\frac{c(D)-1}{2}=\frac{c(K)-1}{2}}$ we have $\displaystyle{u(D)=\frac{c(D)-1}{2}}$. Then by (1) we have that $D$ is a diagram of some $(2,p)$-torus knot as desired. $\Box$

\begin{figure}[htbp]
\begin{center}
\scalebox{0.5}{\includegraphics*{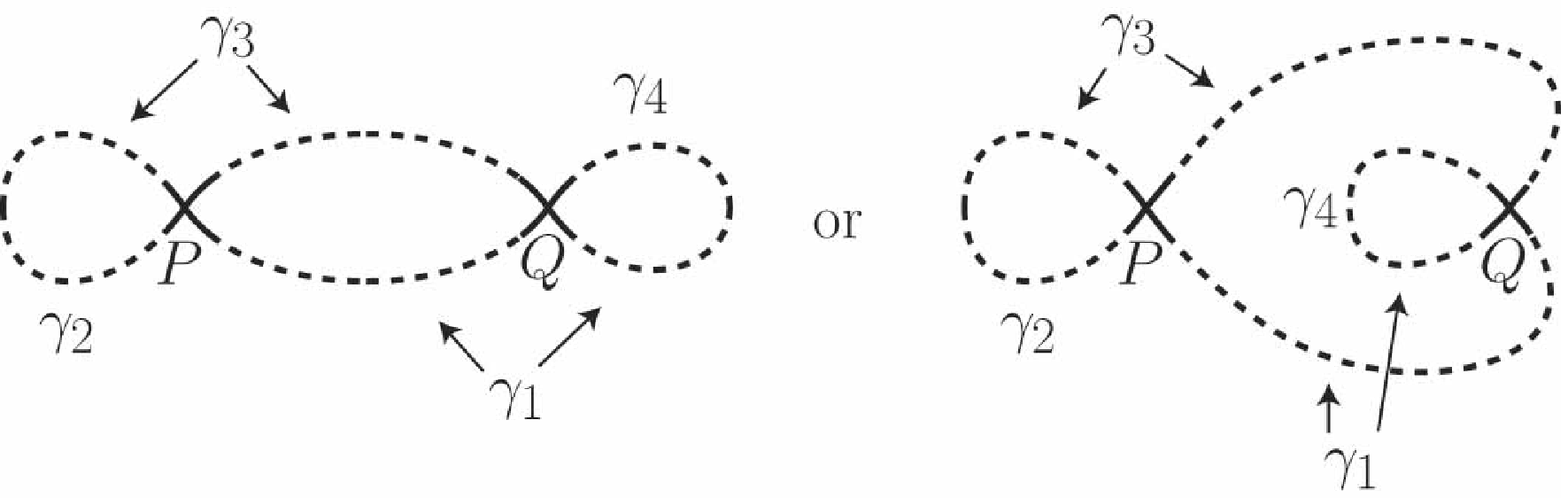}}
\end{center}
\centerline{Figure 3.1}
\end{figure} 

\vskip 3mm

\begin{figure}[htbp]
\begin{center}
\scalebox{0.4}{\includegraphics*{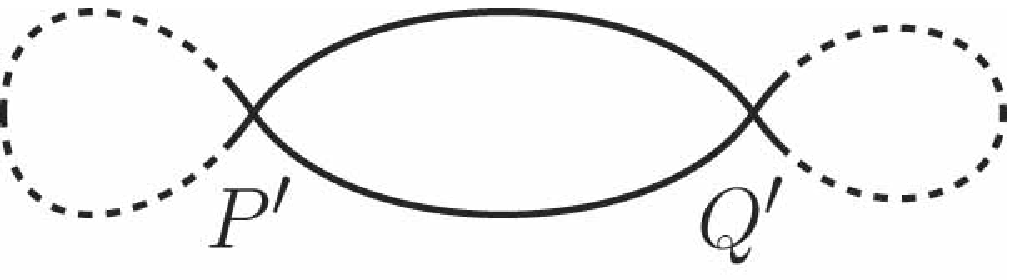}}
\end{center}
\centerline{Figure 3.2}
\end{figure} 

\vskip 3mm

\begin{figure}[htbp]
\begin{center}
\scalebox{0.6}{\includegraphics*{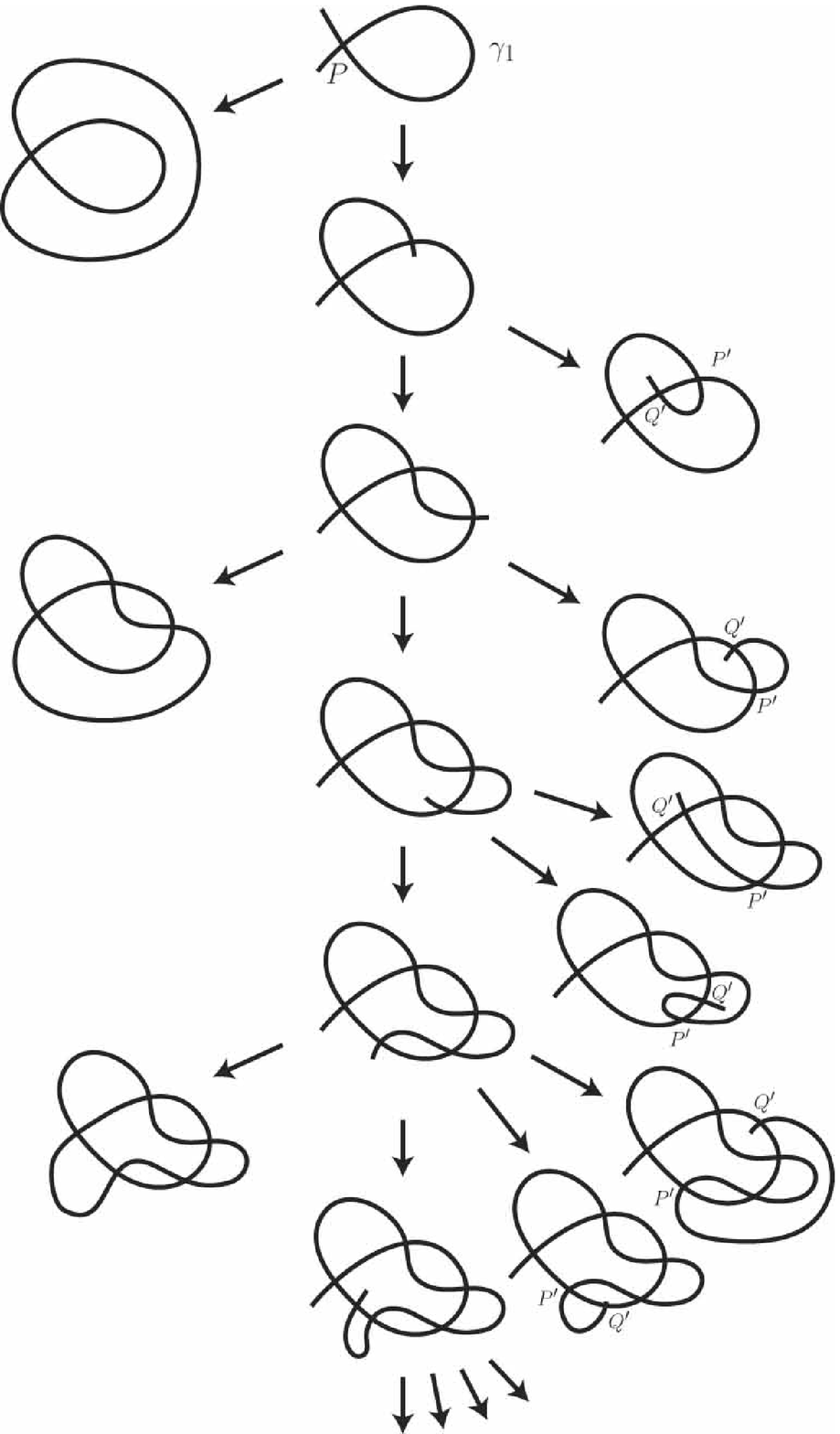}}
\end{center}
\centerline{Figure 3.3}
\end{figure} 

\vskip 3mm

\vskip 3mm 
\noindent{\bf Proof of Theorem 1.5.} Let $D$ be a diagram of a link $L$. Suppose that $D$ has a crossing of the same component. Then as in the proof of Theorem 1.4 we have that $\displaystyle{u(D)\leq\frac{c(D)-1}{2}}$. Therefore we may suppose that every component of $D$ is a simple closed curve on ${\mathbb S}^2$. Suppose that two of them, say $\gamma_1$ and $\gamma_2$ form a non-alternating subdiagram $\gamma_1\cup\gamma_2$. Then by changing at most $\displaystyle{\frac{c(D)-c(\gamma_1\cup\gamma_2)}{2}}$ crossings of $D$ we remove other components and have a link with a diagram $D'$ that consists of $\gamma_1\cup\gamma_2$ and some trivial circles. Then we have $c(D')=c(\gamma_1\cup\gamma_2)$ and $u(D')=u(\gamma_1\cup\gamma_2)$. Since the diagram $\gamma_1\cup\gamma_2$ is non-alternating we have a situation as illustrated in Figure 3.4. Then by changing at most $\displaystyle{\frac{c(\gamma_1\cup\gamma_2)-2}{2}}$ crossings of $\gamma_1\cup\gamma_2$ we have a trivial link of 2-components. The key point here is that we do not need to change two crossings illustrated in Figure 3.4. Therefore we have that $\displaystyle{u(\gamma_1\cup\gamma_2)\leq\frac{c(\gamma_1\cup\gamma_2)-2}{2}}$. After all we have $\displaystyle{u(D)\leq u(D')+\frac{c(D)-c(\gamma_1\cup\gamma_2)}{2}\leq \frac{c(D)-2}{2}}$.

(2) Let $L$ be a link that satisfies the equality $\displaystyle{u(L)=\frac{c(L)}{2}}$. Let $D$ be a minimal crossing diagram of $L$. Since $u(L)\leq u(D)$ and $\displaystyle{\frac{c(D)}{2}=\frac{c(L)}{2}}$ we have $\displaystyle{u(D)=\frac{c(D)}{2}}$. Then by (1) we have that $D$ is a diagram described in Theorem 1.5 as desired. $\Box$

\begin{figure}[htbp]
\begin{center}
\scalebox{0.6}{\includegraphics*{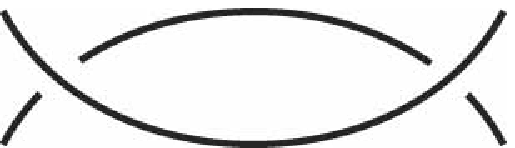}}
\end{center}
\centerline{Figure 3.4}
\end{figure} 

\vskip 3mm

\section*{Acknowledgments}
This work was done during the author was visiting the Department of Mathematics of the George Washington University. He would like to express his hearty gratitude to the Department for its hospitality. In particular the author is grateful to Professor J\'{o}zef Przytycki for his constant guidance and encouragement since 1990. The author is also grateful to Professor Makoto Ozawa for his helpful comments. Finally the author would like to express his gratitude to the referee that his or her comments produced Theorem 1.4 and Theorem 1.5.

\end{document}